\documentstyle[12pt]{article}

\newcommand{\nc}{\newcommand}

\newcommand{\eqdef}{\;\: {\stackrel{ {\rm def} }{=} } \;\:}

\newcommand{\beq}{\begin{equation}}
\newcommand{\eeq}{\end{equation}}
\newcommand{\beqst}{\begin{equation*}}
\newcommand{\eeqst}{\end{equation*}}
\newcommand{\barr}{\begin{array}}
\newcommand{\earr}{\end{array}}
\newcommand{\beqar}{\begin{eqnarray}}
\newcommand{\eeqar}{\end{eqnarray}}
\newtheorem{theorem}{Theorem}[section]
\newtheorem{corollary}[theorem]{Corollary}
\newtheorem{lemma}[theorem]{Lemma}

\newtheorem{definition}[theorem]{Definition}
\newtheorem{remit}[theorem]{Remark}

\newcommand{\CC}{{{\bf  C }}}
\nc{\FF}{ {\Bbb F} }

\newcommand{\calf}{{\mbox{$\cal F$}}}

\newcommand{\call}{{\mbox{$\cal L$}}}
\newcommand{\calm}{{\mbox{$\cal M$}}}

\newcommand{\calo}{{\mbox{$\cal O$}}}

\newcommand{\calr}{{\mbox{$\cal R$}}}
\newcommand{\cals}{{\mbox{$\cal S$}}}

\def\g{\gamma}

\begin{document}

\nc{\Proof}{ \noindent{\bf Proof:} }
\nc{\ctwo}{{F^c}}
\nc{\qtwo}{{F^q}}
\nc{\orb}[1]{{ \calo_{#1} } } 
\nc{\redpar}{{t}}
\nc{\altorpar}{\Lambda}
\nc{\orpar}{{\Lambda}}
\nc{\cone}{Q}
\nc{\proj}{\Pi}
\nc{\push}[1]{\iota_{#1} }
\nc{\diffop}{{\call}}
\nc{\dr}[1]{D_{#1}}
\nc{\Conv}{{ \rm Conv} }
\nc{\cthree}{ {S^c}}
\nc{\liek}{{ \bf {g}}}
\nc{\liets}{\liet^*}
\nc{\lieks}{\liek^*}
\nc{\lietpl}{{ \liet_+}}
\nc{\lietplo}{{ \liet^o_+}}
\nc{\lietspl}{\liets_+}
\nc{\lietsplo}{(\liets_+)^o }
\nc{\modclass}{{\cals^c} }
\nc{\qthree}{{\volsppar{0}{3}}}
\nc{\modquant}{{ \cals^q}}
\nc{\threeorpar}{{ \orpar_1, \orpar_2, \orpar_3 }}
\nc{\threealtorpar}{{ \altorpar_1, \altorpar_2, \altorpar_3 }}
\nc{\threeredpar}{{ \redpar_1, \redpar_2, \redpar_3 }}
\nc{\group}{{G}}
\nc{\torus}{{ T}}
\nc{\domp}{ {D_+} }
\nc{\Del}{\bigtriangleup} 
\nc{\afun}{a}
\nc{\atil}{\tilde{a} } 
\nc{\charac}[1]{\epsilon(#1)}
\nc{\intlat}{\Lambda^I}
\nc{\winv}{{ \frac{1}{|W|} }}
\nc{\conjclass}{ {\rm Cl} } 
\nc{\modsp}{ { \calm} }
\nc{\modsppar}[2]{\modsp_{#1,#2} } 
\nc{\volsppar}[2]{S_{#1,#2} } 
\nc{\volparone}{\volsppar{g}{1}(\Lambda)}     
\nc{\volpar}{\volsppar{g}{n} (\bt) } 
\nc{\modspparreg}[2]{{\modsp^o_{#1,#2} } }
\nc{\upmodsppar}[2]{\calr_{#1,#2} } 
\nc{\upmodspparreg}[2]{\calr^o_{#1,#2} } 
\nc{\sig}[2]{\Sigma^{#1}_{#2} }
\nc{\tc}[1]{[{#1}] }
\nc{\mommap}[1]{\mu^{#1}}
\nc{\mommapo}{{\mu}}
\nc{\modpar}{\modsppar{g}{b} (\Lambda^{(1)}, \dots, \Lambda^{(b)}) } 
\nc{\modparc}{\modsppar{g}{b} (c, \bfl) } 
\nc{\modparone}{\modsppar{g}{1} (\Lambda) } 
\nc{\modparonec}{\modsppar{g}{1} (c,\Lambda) } 
\nc{\extmodsppar}[2]{\tilde{\calm}_{#1,#2} } 
\nc{\extmodpar}{\extmodsppar{g}{n} } 
\nc{\extmodparone}{\extmodsppar{g}{1} } 
\nc{\inprd}[1]{{ (#1) }}
\nc{\inprt}[1]{{ \langle #1 \rangle }}
\nc{\curve}[1]{{ C_{#1} }}
\nc{\momexto}{{J}}
\nc{\momext}[1]{{\momext_{#1} }}
\nc{\alc}{{ D_+} }
\nc{\alco}{{ D^o_+} }
\nc{\Ad}{ {\rm Ad}}
\nc{\bt}{ {\bf { \Lambda}} }
\nc{\bd}{ {c} }
\nc{\mapsing}{{ \phi }}
\nc{\tbunsppar}[3]{{ V_{#1,#2}^{(#3)}  } }
\nc{\tbunspparcirc}[4]{{ V_{#1,#2,#3}^{(#4)}  } }
\nc{\tbunpar}{\tbunsppar{g}{n}{m} (\bt) }
\nc{\tbunparone}{\tbunsppar{g}{1}{m} (\Lambda) }
\nc{\tbunparnn}{\tbunsppar{g}{n}{1} (\bt) }
\nc{\tbunparcirc}{\tbunspparcirc{g}{n}{\alpha}{m} (\bt) }
\nc{\tbunparcircg}{\tbunspparcirc{g}{n}{\gamma_j}{m} (\bt) }
\nc{\loopp}[1]{{ [S_{#1} ]    }}
\nc{\tilal}{\tilde{\alpha} }
\nc{\lietp}{{ \liet^\perp} }
\nc{\nplus}{{ n_+} }
\nc{\dist}{{ f} }
\nc{\qfour}{{ S_{0,4,f } }}
\nc{\fouraltorpar}{{ \altorpar_1, \altorpar_2, \altorpar_3, \altorpar_4 } }
\nc{\waff}{ { W_{\rm aff} }}
\nc{\torpar}{{ h}}
\nc{\tps}{{ P }}
\nc{\qtw}{{ q}}
\nc{\dimg}{{d}}
\nc{\dimt}{{ l}}
\nc{\pant}{ {P}}
\nc{\triup}{ { \bigtriangleup} } 
\nc{\free}{ { \Bbb F } } 
\nc{\univ}{ {\Bbb U} } 
\nc{\flagm}{ {\rm Fl} }
\nc{\hattau}{ { \hat{\tau} } }
\nc{\hu}[1]{ {\hat{u}_{#1} } }
\nc{\modonetwo}[1]{ \calm_{1,1}^{{SU(2)}} (#1) }
\nc{\tortwo}{ { T_2 } }
\nc{\tbunonetwo}[1]{ V_{1,1}^{{SU(2) }} (#1) }
\nc{\subv}{ { Z} } 
\nc{\tork}{ { T_k} }
\nc{\tilzed} {\tilde{Z} } 
\nc{\zed}{ {Z} }
\nc{\inclzed} {i_{\zed} }
\nc{\divis}[3]{D_{#1, #2}(#3)}
\nc{\tildediv}[3]{\tilde{D}_{#1, #2}(#3)}
\nc{\lamax}{ \bigwedge^{\rm max} }
\nc{\nminus}{{N^-}}
\nc{\omkd}{{ \omega_{n,d} }}
\nc{\secc}[2]{s_{#1}^{(#2)} }
\nc{\partlam}{\partial_\Lambda}
\nc{\partlamfun}[1]{\partial_{\Lambda_{#1}}}

\nc{\he}[1]{{\hat{u}_{#1} } }
\nc{\uroot}[1]{ {u_{#1}}  }
\newcommand{\labell}{\label}
\input amssym.def
\input amssym.tex
\newcommand{\renorm}{{ \setcounter{equation}{0} }}
\renewcommand{\theequation}{\thesection.\arabic{equation}}
\nc{\mnd}{{ M(n,d) } }
\nc{\weightl}{\Lambda^W}
\nc{\liet}{ {\bf t} } 
\nc{\lieg}{ { \bf g} } 
\nc{\ct}{ { \tilde{c} } }
\nc{\tilc }{ { \tilde{c} } }
\nc{\diag}{{ \rm diag} } 
\nc{\Res}{{\rm Res} } 
\nc{\res}{{\rm Res} } 
\nc{\lb}[1]{ {l_{#1} } }
\nc{\yy}[1]{Y_{#1}}
\nc{\liner}[1]{L_{#1} }
\nc{\lambdr}[1]{\Lambda_{#1} }
\nc{\linestd}{\liner{(\lb{1}, \dots, \lb{n-2} )} }
\nc{\lambstd}{\lambdr{(\lb{1}, \dots, \lb{n-2} )} }
\nc{\expsum}[1]{ { (e^{{#1} } - 1 ) } }
\nc{\itwopi}{ { 2 \pi i }}
\nc{\indset}{{(l_1, \dots, l_{n-1})}  }
\nc{\indsettwo}{{(l_1, \dots, \dots, l_{n-2})}  }
\nc{\bfl}{{ \bf \Lambda}}
\nc{\dgk}{{ D_{n,d}(g,k) }} 
\nc{\vgkl}{{ D_{n,d}(g,k, \Lambda) }} 
\nc{\vgklmult}{{ D_{n,d}(g,k, \bfl) }} 
\nc{\dgklmult}{{ V_{n,d}(g,k, \bfl) }} 
\nc{\vgk}{{ V_{n,d}(g,k) }} 
\nc{\dgkl}{{ V_{n,d}(g,k,\Lambda) }} 
\nc{\lineb}{L}
\nc{\resid}[1]{\res_{#1 = 0 } }
\nc{\residone}[1]{\res_{#1 = 1 } }
\nc{\kmod}{r}
\nc{\gmax}{{\g_{\rm max} }}
\nc{\constq}{ {\frac{1}{n}  }}
\nc{\sol}{ { S_{0 \mu} }}
\nc{\ch}{ { \rm ch} } 
\nc{\td}{ { \rm Td} } 
\nc{\abk}{\kappa}
\nc{\tg}{{\tilde{\gamma} } }
\nc{\inpr}[1]{ \langle #1 \rangle }
\nc{\ib}{\item{$\bullet$} }
\newcommand{\sg}{\Sigma^g} 
\title{The residue formula and the Tolman-Weitsman theorem}
\author{Lisa C. Jeffrey \\
Mathematics Department \\
University of Toronto\\
Toronto, Ontario M5S 3G3, Canada
\thanks{This material is based on work
supported by 
 grants from NSERC and the Alfred P. Sloan Foundation.
 \hspace*{\fill} MSC subject classification: 58F05} }
\date{ May 2002} 
\maketitle
\begin{abstract}

We give a simple direct proof (for the case of Hamiltonian
circle actions with isolated fixed points) that Tolman and Weitsman's
description of the kernel of the Kirwan map 
in \cite{TW}  (in other words the sum
of those
equivariant cohomology classes vanishing on one side of a collection 
of hyperplanes) is equivalent to the characterization of this
kernel  given by
the residue theorem \cite{JK}. 

\end{abstract}

\newcommand{\hht}{{H^*_T}}
\newcommand{\mred}{{M_{c}}}

\nc{\fplus}{\calf_+}
\nc{\fminus}{\calf_-}
\nc{\almin}{{\alpha^-}}
\nc{\alplus}{{\alpha^+}}

\section{Introduction}

Let $M$ be a symplectic manifold of dimension $2n$ equipped with the 
Hamiltonian action  of a torus $T$. The $T$-equivariant cohomology
of $M$ is $H^*_T(M)$. This is a module over $\hht({\rm pt}).$ 
The moment map for the torus action is denoted by
$ \mu: M \to \liets$.
We denote by $\mred$ the reduced space at a value $c$ of the moment 
map: 
$$ \mred \eqdef \mu^{-1}(c)/T. $$

The Kirwan map is the  map
$$\kappa: \hht(M) \to H^*(\mred) $$
induced by the restriction map
$$\hht(M) \to \hht\left (\mu^{-1}(c)\right ). $$
Provided $c$ is a regular value of $\mu$ (which we shall assume)
we have that 
$\hht\left (\mu^{-1}(c)\right )  \cong H^*(\mred)$.\footnote{In
this paper all cohomology groups are with complex coefficients.}

Kirwan proved in \cite{Kir} that the Kirwan map is 
surjective. Thus to find the cohomology ring of 
$\mred$ it suffices to find the 
kernel $K$ of $\kappa$:
$$H^*(\mred) \cong \hht(M)/K. $$ 
We denote by $\calf$ the set of components of the fixed point set of $T$.
We define a distinguished  subset $\fplus$ of $\calf$: in the case when 
$T = S^1$ it is simply 
$$ \fplus = \{ F \in \calf | \mu(F) > c\}. $$
In the case $T = S^1$ we also introduce 
$$ \fminus = \{ F \in \calf | \mu(F) < c \}. $$
Because 
$c$ is a regular value of $\mu$, there are 
no $F$ with $\mu(F) = 0 $.

In \cite{JK} the following is proved:

\begin{theorem} \label{t:jkres}
Let $\eta \in \hht(M)$. Then
$$\kappa (\eta)[\mred] = \sum_{F \in \calf_+} {\rm Res} \left (
\frac{\eta}{e_F}\right ) [F]. 
$$  
\end{theorem}
Here $e_F$ is the equivariant Euler class of the 
normal bundle to $F$, and 
${\rm Res} $ is an iterated residue.\footnote{For more details
see Section 8 of \cite{JK}: see also \cite{locq} and \cite{GK}.}
In the case 
when $T$ has rank one it is simply ${\rm Res}_{X = 0 }$ where
the variable 
$X $ is the generator of $\hht({\rm pt})$.

Since $\kappa$ is a ring homomorphism, we have from Theorem \ref{t:jkres}
\begin{corollary}
\beq \label{eq:1}
\kappa(\eta) \kappa( \zeta) [\mred] = 
\sum_{F \in \calf_+} {\rm Res} \left (\frac{\eta \zeta}{e_F}\right )[F] . \eeq 
\end{corollary}
Assuming $c$ is a regular value of $\mu$, 
$\mred$ has at worst orbifold singularities and so satisfies
Poincar\'e duality. 
Thus according
to Theorem 
\ref{t:jkres}, an element $\eta \in \hht(M)$ is in $K$ if and only 
if 
$$\sum_{F \in \calf_+} {\rm Res} \left(\frac{\eta \zeta}{e_F}
\right ) [F]  = 0 $$
for all $\zeta \in \hht(M)$. 

In \cite{TW} Tolman and Weitsman give an alternative characterization 
of the kernel $K$:

\begin{theorem} \label{t:tw} 
The kernel $K$ of the Kirwan map is the sum
$$K = \sum_{\xi \in \liet} K_\xi$$
where $K_\xi$ is  the   set of cohomology classes vanishing on 
 all components $F$ of the 
fixed point set for which
$(\mu(F)-c, \xi) > 0 $. 
\end{theorem}

In the special case $T = S^1$  the Tolman-Weitsman theorem
reads as follows:

\begin{corollary}\label{c:tw1}
Let $T = S^1$. Then
the kernel $K$ of the Kirwan map is the sum
$$K = K_+ \oplus K_- $$
where $$K_{\pm}  = \{ \eta \in \hht(M)~|~ 
\eta[F] = 0 \mbox{ for all}~ F\in \calf~
\mbox{for which}~ \pm \mu(F) > c\}.  $$
\end{corollary}

Tolman and Weitsman prove their theorem by 
establishing directly that 
$\sum_\xi K_\xi$ is the kernel of the Kirwan 
map $\kappa$, and that all elements of this kernel take this form.

In this paper  we give a direct proof in the case $T = S^1$ that
the classes satisfying
the condition that  the  residue on the
right hand side of (\ref{eq:1}) is
zero  are contained in the Tolman-Weitsman kernel
$\oplus_\xi K_\xi$.
The reverse inclusion is clear from the 
definition of the residue: 
if a class $\eta \in \hht(M)$ satisfies $\eta|_F = 0 $ for 
$F \in \calf_+$, then clearly 
\beq \label{e:res} 
\sum_{F \in \calf_+} {\rm Res}_{X = 0 } 
\left ( \frac{\eta \zeta}{e_F}\right )  [F] = 0 \eeq
 for all 
$\zeta \in \hht(M)$. Similarly if $\eta$ satisfies 
$\eta|_F = 0 $ for all $F \in \calf_-$ then the residue 
in (\ref{e:res}) is equal to zero.
Thus we have shown that the description of the kernel of the 
Kirwan map given by the residue theorem is equivalent to the 
description given by Tolman and Weitsman. We expect that
our methods could be enlarged to treat the case of general 
torus actions.

\section{Results from Morse theory}

We begin by stating a number of fundamental
results related to equivariant Morse theory and the 
Thom-Gysin map. These results
are found in the work of Atiyah and Bott \cite{AB}
and Kirwan \cite{Kir}: the formulation we  present
appears in the work of Goldin \cite{Gold1}.
We introduce a partial ordering $<$ on the components
of the critical set of a Morse function $f$:
denote these components by $F_i$. The 
ordering is given by the value of 
the 
Morse function $f$: 
in other words $i > j $ if and only if $f(F_i) > f(F_j)$. 

In our situation the most important 
Morse function $f$ is  $-\mu_\xi$,
where $\mu$ is the moment map and $\mu_\xi = (\mu,\xi)$ is the 
component of the moment map in the direction $\xi$ for nonzero values $\xi \in \liet$.  
This is an equivariantly  perfect Morse function
(as proved by Atiyah-Bott and Kirwan \cite{AB,Kir}).

We denote by  ${\rm Ind}(F) $ 
the Morse index, in other words the number of negative eigenvalues
of the Hessian of $f$ on the critical set $F$. 

In this paper we assume for simplicity that all components $F$ of the 
critical set of $f$ are isolated (so $f$ is a Morse function rather than 
a Morse-Bott function). 

\begin{theorem}\label{t:g2.4}Suppose 
$\eta \in \hht(M)$ restricts to $0$ on all 
$F_i$ for which $i < j$. Then $\eta|_{F_j}$ is some
multiple of $e \left (\nu_f^- F_j \right ) $, the equivariant
Euler class of the negative normal
bundle of $F_j$. 
\end{theorem}

\begin{remit} \label{r:deg} An important consequence
of 
Theorem \ref{t:g2.4} is that if $\eta$ restricts
to $0$ on all $F_i $ with $i < j$ then 
either $\eta|_{F_j} = 0 $ or the degree of 
$\eta$ is greater than or equal to the
degree of $e(\nu^-(F_j))$.
\end{remit}


\begin{definition}
We say that an equivariant cohomology class
$\zeta \in H^*_{S^1}(M) $ is supported on a subset
$\cals$ of the critical set of $f$ if $\zeta|_G = 0 $ for all $G $ in the 
critical set for which
$G \notin \cals$. 
\end{definition}

\begin{definition} \label{flowup}
Let $F$ be a component of the critical set of $f$. The 
 equivariant stable manifold  of $F$ under $f$
 is the set of points $x \in M$ for which
there is a trajectory  which
converges to $F$ under the flow of $- {\rm grad} ( f)$.
\end{definition}

\begin{definition} \label{extflowup}
Let $F$ be as in Definition \ref{flowup}. The equivariant extended 
stable manifold
of $F$ under $f$ is the set of points $x \in M$ which have a sequence of 
trajectories of the flow of $- {\rm grad} (f)$ (passing through 
a number of components $F_1, \dots, F_l$ of the 
critical set of $f$) starting at  $x$ and 
converging finally to  $F$.
\end{definition}

\begin{definition}
Let $F$ be as in Definition \ref{flowup}. The unstable
manifold and equivariant
extended unstable manifold of $F$ are specified as in Definitions
\ref{flowup} and  \ref{extflowup}
but using the gradient flow of ${\rm grad} (f)$ rather than that of 
$- {\rm grad} (f)$.
\end{definition}

\begin{theorem}\label{t:g2.5}

Suppose $F$ is a connected component of 
the critical set of a Morse function $f$. 
Then there  is a class $\almin(F)$ with the following properties:
\begin{enumerate}
\item $\almin(F)|_G =  0$ if $G$ is a component of the critical 
set of $f$ which is not in the equivariant extended stable manifold of $F$ 
\item $\almin(F)|_F = e(\nu_f^- F) $ where $e(\nu_f^- F)$ is 
the equivariant Euler class of the negative normal bundle
(defined by $f$) of $F$.
Thus the degree of $\almin(F)$ is {\rm Ind}(F). 
\end{enumerate}

In the same way there is a class $\alplus(F)$ 
such that
\begin{enumerate}
\item $\alplus(F)|_G =  0$ if $G$ is a component of the critical 
set of $f$ which is not in the  equivariant extended unstable manifold
 of $F$ 
\item $\alplus(F)|_F = e(\nu_f^+ F) $ where $e(\nu_f^+ F)$ is 
the equivariant Euler class of the positive normal bundle
(defined by $f$) of $F$
\end{enumerate}
The degree of $\alplus(F)$ is $2n- {\rm Ind}(F).  $
\end{theorem}

\section{The $S^1$ case}
  
We now assume that the torus $T$ is the circle group $S^1$, and also that
the fixed points of the action of $T$ are isolated. 
Suppose we are forming the  reduced space at a regular
value $c$ of the moment map. 
The purpose of this paper is to give a direct proof of the following 
result.

\begin{theorem} \label{t:sone}
Let $\eta \in \hht(M)$ satisfy the hypothesis that 
$$\sum_{F \in \calf_+}  {\rm Res}_{X = 0 } \left ( 
\frac{\eta \zeta}{e_F}\right ) [F]  = 0 $$
for all $\zeta \in \hht(M)$. Then 
$ \eta \in K_+ \oplus K_-$. 
\end{theorem}

\nc{\almaxminus}{\alpha_{\rm max}^-}
\nc{\kplus}{{ K_+}}
\nc{\kminus}{{ K_-}}
\nc{\gp}{{S^1}}

\noindent{\bf Proof of Theorem \ref{t:sone}:}
Let $\eta \in H^d_{\gp} (M)$, in other words $d $ is the degree of 
$\eta$.

\begin{lemma} We can write
\beq \label{e:decomp}\eta = \eta_+ + \eta_- \eeq
where 
\beq  \label{e:decompdefmin}
\eta_- = \sum_{F \in \fminus, {\rm Ind}(F) \le d}  c_F X^{n_F}
\almin (F) \eeq
while
\beq \label{e:decompdefplus}
 \eta_+ = \sum_{F \in \fplus, {\rm Ind}(F) \le d}  c_F X^{n_F} \almin(F) \eeq
and 
\beq \label{nfdef} 2 n_F = d - {\rm Ind}(F). \eeq 
Here, the $c_F \in \CC $ are constants. 
\end{lemma}

\Proof This is straightforward. We simply have to identify the 
coefficients of the restrictions to components $F$ of the 
fixed point set. We find that we can adjust the coefficients in such 
a way that both sides have the same restriction to all components
of the fixed point set, which by Kirwan's injectivity theorem \cite{Kir}
(as well as Remark \ref{r:deg})
guarantees
that the two are equal in equivariant cohomology.

We define the constants $c_F$ by induction using $\almin(F)|_F = 
e(\nu_f^- F)$.
 Adding
$c_F X^{n_F} \almin(F)$ will alter only the values of the restrictions
to those elements $G\in  \calf$ for which $\mu(G) \ge \mu(F)$.  \hfill
$\square$

\begin{lemma} \label{l3.7} We have that 
$$ \eta_+ \in \kminus. $$
\end{lemma}
\Proof For $F \in \fplus$, $\almin(F)$ is supported on 
$\fplus$ (since if $G$ is a component of $\calf$ for which
$\alpha(F)|_G \ne 0 $, $G$ is in the equivariant extended 
stable manifold
of $F$ so $\mu(G) \ge \mu(F)$). \hfill $\square$

\begin{lemma} \label{l3} 
We may add terms
$$b_F X^{n_F} \almin(F) $$
for $F \in \fplus$, in such a way that 
$$\eta_- - \sum_{F \in \fplus} b_F X^{ {n_F}} \almin(F) $$
restricts to zero on $G$ for all $G \in \fplus$ for which ${\rm Ind} (G) 
\le d$. Here, $n_F$ is as in (\ref{nfdef}). 
\end{lemma}
\Proof We choose the coefficients $b_F$ by induction on 
the value of $\mu(F)$), first choosing $b_F$ for the lowest value of 
$\mu(F)$  for $F \in \fplus$. This will alter only the value of the restrictions
to those $G$ for which $\mu(G) \ge \mu(F)$. 
This process can be used to make the values of $\eta_- (G) = 0 $
for all $G$ for which $\mu(G) > c $ and ${\rm Ind}(G) \le d$. 
\hfill $\square$

\nc{\alphaplus}{{\alpha_+}}
\nc{\alphaminus}{{\alpha_-}}

As described in the proof of  Lemma \ref{l3.7}, 
for $G \in \fplus$, $\alpha^-(G)$ is  supported on $\fplus$. So we have
the following corollary:
\begin{corollary} \label{c:6}
Suppose $G \in \fplus$. Then
we see that $ \alphaminus(G) \in \kminus$ and thus
\beq 
\label{e:fplusvan}{\rm Res}_{X = 0 }  \sum_{G \in \fminus} 
\frac{\alphaminus (G)|_G ~
\zeta|_G}{e_G}
 =0 \eeq
for any  class $\zeta \in \hht(M)$.
Therefore also 
\beq \label{e:fminvan}{\rm Res}_{X = 0 } 
\sum_{G \in \fplus} \frac{\alphaminus(G)|_G ~  \zeta|_G}{e_G} = 0, \eeq
since the sum of residues for $F \in \fminus$ is equal to minus
 the sum of residues for $F \in \fplus$ for any element of $\hht(M)$.
Because 
$$
\sum_{G \in \fplus \cup \fminus} \frac{\alphaminus(G)|_G ~ \zeta|_G}{e_G} $$
 is a smooth function of $X$
(by the localization theorem in equivariant
cohomology \cite{AB2},\cite{BV}), it 
does not have any residue at $X = 0 $,
so (\ref{e:fplusvan}) implies (\ref{e:fminvan}).
\end{corollary}
\begin{lemma} \label{l7}  
Write 
$$ \eta_-|_F = a_F X^{d/2} $$
for some $a_F \in \CC$. 
If $\eta_-|_F = 0 $ for all $F\in \fplus$ for which 
${\rm Ind}(F) \le d$, then 
also $\eta_-|_F = 0 $ for all $F \in \fplus$. In other
words $\eta_- \in K_+. $  
\end{lemma}
\begin{remit} 
The hypothesis of Lemma \ref{l7} can be achieved because of Lemma \ref{l3}. 
\end{remit}
\noindent{\bf Proof of Lemma \ref{l7}:}

By Corollary \ref{c:6}, 
 we may assume 
$${\rm Res}_{X = 0 } \sum_{G \in \fplus} \frac{\eta_-|_G ~\zeta|_G}{e_G} = 0 $$
for all classes $\zeta$. We put a partial order on 
 the elements $F \in \fplus$ according to the 
values of $\mu(F)$, as described
above.
We choose 
\beq \label{zetdef}\zeta  = X^{m_F} \alphaplus(F) \eeq
for 
\beq \label{e:bound1} 2n - {\rm Ind} (F) + d \le 2n - 2 \eeq	 
and an appropriate nonnegative integer $m_F $ 
(in other words
\beq \label{e:bound2} {\rm Ind}(F) \ge d + 2:  \eeq
this hypothesis is needed in order
to obtain terms with nonzero residue in the sum). 
We choose $m_F$ so that
${\rm deg} (\zeta \eta) = 2n-2,$ 
in other words
\beq  \label{e:3.8}  2 m_F  +2 + d  =  {\rm Ind}(F). \eeq
We allow $\zeta$ to be determined by (\ref{zetdef}) and
(\ref{e:3.8}) for all 
$F$ for which
${\rm Ind}(F) \ge d+2.$ 
Notice that because the fixed point set consists of isolated points,
Kirwan's injectivity theorem implies that all $\eta \in H^*_T(M)$
have even degree. Hence there are no components $F$ 
for which $\eta|_F \ne 0 $ and 
${\rm Ind}(F) = d + 1$ (using (\ref{e:3.8})).

For all $G \in \fplus$ we define $b_G(F) \in \CC$ by 
\beq \label{star} \alphaplus(F)|_G = b_G(F) X^{n-({\rm Ind}(F)/2)}. \eeq
Note that $b_G(G) \ne 0 $ since $\alpha_+(G)|_G = e_G^+(G) \ne 0$. 
We also write
\beq \label{e:3.9a} e_F|_F = \epsilon_F X^n \eeq 
for some nonzero $\epsilon_F \in \CC$. 
We define $a_G \in \CC$ by
\beq \label{dagger} \eta_-|_G = a_G X^{d/2}. \eeq
The equation 
$${\rm Res}_{X = 0 } \sum_{G \in \fplus} 
\frac{\eta_-|_G ~ \zeta|_G} {e_G}  = 0 $$
then becomes (combining (\ref{zetdef}), (\ref{star}),
(\ref{e:3.9a}) and (\ref{dagger}))
\beq \label{e:master}
\sum_{G \in \fplus} \frac{a_G~ b_G(F)}{\epsilon_G}  = 0 . \eeq
For the  value of  $F $  for which 
$\mu(F) $ is maximal, the equation (\ref{e:master}) is simply
\beq \sum_{G \in \fplus} \frac{a_G}{\epsilon_G}  = 0, \eeq
since in this case the class $\alphaplus (F)$ has degree zero 
and is represented in the Cartan model by 
a constant in 
$\CC$, so the $b_G(F)$ all take the same value.

Notice that $b_G(F) = 0 $ if $G$ is not in the equivariant extended 
unstable manifold of $\mu$ for $F$ (for example, if $\mu(G) > \mu(F)$ then 
$b_G(F) = 0 $).
Thus $b_G(F) \ne 0 $ only for $G$ in the equivariant extended unstable
manifold
of $F$. 
By Lemma \ref{l3}, we  can also assume  
that $a_G = 0 $ whenever ${\rm Ind}(G) \le d$.
Also if ${\rm Ind}(G) = d+ 1$ then $a_G = 0 $ since ${\rm Ind}(G)$  must
be odd so $\eta|_G = 0 $ (by {\ref{e:3.8}). 
Thus the residue formula implies that 
(\ref{e:master}) holds for all $F$ for which
${\rm Ind}(F) \ge d+2$.   By (\ref{e:bound1})
and (\ref{e:bound2}), we see that
(\ref{e:master}) becomes
the equation 
\beq \label{e:master2}
\sum_{G \in \fplus, ~{\rm Ind}(G) \ge d+2 }
 \frac{a_G~ b_G(F)}{\epsilon_G}  = 0  \eeq
which is valid whenever ${\rm Ind} (F) \ge d + 2$.

To deduce from this that $a_G = 0 $ for all $G \in \fplus$, we think
of (\ref{e:master}) as a matrix equation involving an $N \times N$  matrix
$B$ 
(where $N$ is the number of elements $F \in \fplus$
with ${\rm Ind}(F) \ge d+2$): the matrix is 
$$ B_{FG} = b_G(F) $$
for all $F$ and $G$ with ${\rm Ind}(F) \ge d+2$ and 
${\rm Ind}(G) \ge d+2$. 
(Here the $F$ index the rows and $G$ index the columns of the matrix.)
We wish to solve (\ref{e:master}) for the vector $(a_1, \dots, a_N)$.
We note that $b_G(F) = 0 $ if $\mu(G) > \mu(F) $ or $\mu(G) = \mu(F)$ and
$G$ is not in the equivariant extended  unstable
manifold of $F$. Thus 
the matrix $B$ is upper triangular; furthermore its entries on the 
diagonal are nonzero. 
Thus the matrix is invertible, so the only solution to the
equations (\ref{e:master}) is 
$a_G = 0 $ for all $G$. This completes the proof 
of Lemma \ref{l7}. \hfill $\square$

\noindent{\bf Proof of Theorem \ref{t:sone}:} We now have completed
the proof of the Theorem, since we have shown that $\eta $ can be
written 
$$\eta = \eta_+ + \eta_-$$
where it is possible to find $\eta_- \in K_+$  and  $\eta_+ \in K_-$.
\hfill $\square$


\begin{thebibliography}{99}
\bibitem{AB} M.F. Atiyah and R. Bott, The Yang-Mills equations
over Riemann surfaces, {\em Phil. Trans. Roy. Soc. Lond.}
{\bf A308} (1982) 523-615.
\bibitem{AB2} M.F. Atiyah and R. Bott, The moment map and equivariant
cohomology, {\em Topology} {\bf 23} (1984) 1-28.
\bibitem{BV} N. Berline and M. Vergne, Z\'eros
d'un champ de vecteurs et classes caract\'eristiques
\'equivariantes, {\em Duke Math. J.} {\bf 50} (1983) 539-549.
\bibitem{Gold1} R. Goldin, {An effective algorithm for the cohomology
ring of symplectic reductions}, preprint math.SG/0110022 (2001);
{\em GAFA}, to appear.
\bibitem{GK} V. Guillemin and J. Kalkman,
The Jeffrey-Kirwan localization theorem and residue operations 
in equivariant cohomology. {\em J. Reine Angew. Math.} {\bf 470} (1996), 
123--142. 
\bibitem{JK} L. Jeffrey and  F. Kirwan,  Localization for nonabelian
group actions, {\em Topology} {\bf 34} (1995) 291-327.
\bibitem{locq} L. Jeffrey and F. Kirwan, 
Localization and the quantization conjecture. 
{\em Topology}  {\bf 36} (1997), 647-693.
\bibitem{Kir} F. Kirwan, {\em Cohomology of Quotients in Symplectic
and Algebraic Geometry}, Princeton University Press, 1984.
\bibitem{TW} S. Tolman and J. Weitsman, {The cohomology ring of 
abelian symplectic quotients}, preprint math.DG/9807173.
\end{thebibliography}
\end{document}